# A THEORETICAL FRAMEWORK FOR THE PRICING OF CONTINGENT CLAIMS IN THE PRESENCE OF MODEL UNCERTAINTY

By Laurent Denis and Claude Martini

*Université d'Evry-Val-d'Essone and INRIA and Zeliade Systems*

The aim of this work is to evaluate the cheapest superreplication price of a general (possibly path-dependent) European contingent claim in a context where the model is uncertain. This setting is a generalization of the uncertain volatility model (UVM) introduced in by Avellaneda, Levy and Paras. The uncertainty is specified by a family of martingale probability measures which may not be dominated. We obtain a partial characterization result and a full characterization which extends Avellaneda, Levy and Paras results in the UVM case.

**1. Introduction.** Our purpose is to set a framework for dealing with model uncertainty in mathematical finance and to handle the pricing of contingent claims in this context.

Let $(S_t)_{t \in [0,T]}$ be a real-valued process which stands for an asset price. Usually, it is assumed that the set $\mathbf{P}_m$ of the equivalent probabilities under which $S$ is a martingale is not empty. This is a sufficient condition to preclude pure gambling strategies that never fail and win with a positive probability.

In the classical example of the Black–Scholes model, $S$ solves the linear stochastic differential equation $dS_t = S_t \, dB_t$, where

$$B_t = \sigma W_t + \nu t,$$

and $W$ is a Brownian motion, $\nu$ is the drift and $\sigma$ is the so-called volatility parameter.

Consider now some bounded random variable $f$ which represents the payoff of a European contingent claim on $S$. The cheapest riskless superreplication price of the claim can be defined as

$$\Lambda(f) = \inf\left\{ a \in \mathbb{R} : \exists h \text{ such that } a + \int_0^T h_t \, dB_t \geq f \ P\text{-a.s.}, \right\}$$

where $h$ is an integrant for which the stochastic integral as a process is uniformly bounded from below (usually the integrator is $S$, but it does not matter because $h_t \, dS_t = h_t S_t \, dB_t$). An important achievement of the theory is the duality formula

$$\Lambda(f) = \sup_{Q \in \mathbf{P}_m} E_Q f,$$

which was obtained first by Delbaen [6] for continuous semimartingales; El Karoui and Quenez [11] derived it from the optional decomposition and their approach was extended by Kramkov [17] to locally bounded processes. The local boundedness assumption eventually was removed by Föllmer and Kabanov [13].

So, for a completely specified model, the superreplication issue is settled. The theory is less satisfactory in the practically important case where the model is not completely specified. The first attempts to attack the problem in such a context were undertaken by Avellaneda, Levy and Paras [2] and Lyons [18], who introduced an uncertainty in volatility (uncertain volatility model, abbreviated UVM): the volatility process is not known, but is assumed to lie in a fixed interval $[\underline{\sigma}, \overline{\sigma}]$. The authors obtained a generalization of the duality formula by stochastic control techniques in the case of European options with payoffs that depend only on the terminal value $B_T$. The discrete-time case has been studied recently in [10].

A major difficulty, which is already present in the UVM case, is that one is faced with a family of measures which are, in general, mutually singular and nondominated (i.e., not absolutely continuous with respect to a single probability measure).

The purpose of this paper is as follows:

- To provide a coherent framework on which one can set the superreplication problem for European contingent claims (including path-dependent ones), which encompasses the case of the UVM model.
- To find a characterization of the cheapest superreplication price.

We work with continuous one-dimensional processes to insulate the technicalities which are specific to our uncertainty framework. In Section 2 we formulate the superreplication problem in the presence of model uncertainty. The main idea is to use capacity theory to define a refined stochastic integral. Then we give our main characterization result. This result is proved first in the case where $B$ is bounded (Section 4) and then proved in the unbounded case (Section 5). In Section 6 we apply our results to a slightly generalized UVM model and obtain an extension of [2] and [18] to path-dependent European contingent claims. In the Appendix we collect necessary facts from capacity theory.

## 2. The superreplication problem in the presence of model uncertainty.

2.1. *A mathematical framework for model uncertainty.* Let $\Omega = C([0,T], \mathbb{R})$—the space of scalar continuous functions $B = (B_t)$ and $B_0 = 0$—be endowed by the uniform norm and let $\mathcal{B}$ be its Borel $\sigma$-field. We denote by $\mathcal{F}_t$ the canonical filtration; $E_P$ stands for the expectation under $P$. A probability $P$ on $(\Omega, \mathcal{B})$ is a *martingale measure* if the coordinate process $B$ is a martingale with respect to $\mathcal{F}_t$ under $P$; we denote by $\mathbf{P}_m$ the set of all such measures. In addition, $\langle B \rangle^P$ stands for the quadratic variation of $B$ under $P$ (it is defined up to a $P$-null set). Fix a nonzero measure $\overline{\mu}$ on $[0,T]$ with continuous distribution function which also will be denoted by $\overline{\mu}$ (i.e., $\overline{\mu}_t = \overline{\mu}([0,t])$). We consider a subset $\mathbf{P} \subseteq \mathbf{P}_m$ of martingale measures that satisfy the following standing.

HYPOTHESIS $H(\overline{\mu})$. For each $P \in \mathbf{P}$, the process $\overline{\mu} - \langle B \rangle^P$ is increasing (up to a $P$-null set). We denote this relationship by

$$d\langle B \rangle^P_t \le d\overline{\mu}_t.$$

This hypothesis ensures that $\mathbf{P}$ is relatively weakly compact, the property believed to be a minimal one for future development (see the Appendix). For certain results we need that $\overline{\mu}$ is Hölder continuous, that is, for some positive constants $C$ and $\alpha$,

(1) $$\overline{\mu}_t - \overline{\mu}_s \le C|t-s|^\alpha, \qquad s,t \in [0,T], s \le t.$$

We shall use also an assumption that involves a second nonzero measure $\underline{\mu}$:

HYPOTHESIS $H(\underline{\mu}, \overline{\mu})$. For each $P \in \mathbf{P}$,

$$d\underline{\mu}_t \le d\langle B \rangle^P_t \le d\overline{\mu}_t.$$

Introducing the above conditions, we have in mind the Black–Scholes model with volatility $\sigma \in [\underline{\sigma}, \overline{\sigma}]$. In such a case, $d\underline{\mu} = \underline{\sigma}^2 \, dt$ and $d\overline{\mu} = \overline{\sigma}^2 \, dt$.

2.1.1. *Definition of the capacity.* For each $f \in C_b(\Omega)$—the set of bounded continuous functions on $\Omega$—we put

$$c(f) = \sup\{\|f\|_{L^2(\Omega, P)} : P \in \mathbf{P}\}.$$

The convex positive homogeneous function $c$, *capacity*, can be extended naturally to all functions on $\Omega$; by definition, $c(A) = c(I_A)$ (for details, see the Appendix). We use the standard capacity-related vocabulary: A set $A$ is *polar* if $c(A) = 0$ [thus, if $A$ is measurable, then $P(A) = 0$ for all $P \in \mathbf{P}$]; a property holds "quasi-surely" (q.s.) if it holds outside a polar set.

A mapping $f$ on $\Omega$ with values in a topological space is *quasi-continuous* if $\forall \varepsilon > 0$, there exists an open set $O$ with $c(O) < \varepsilon$ such that $f|_{O^c}$ is continuous.

REMARK. For applications of capacities that arise from a set of probabilities, we refer the reader to the theory of *risk measures*; see, for example, [1, 7, 14].

2.1.2. *The space of quasi-continuous functions.* We denote by $\mathcal{L}$ the topological completion of $C_b(\Omega)$ with respect to the semi-norm $c$ and denote by $L$ the quotient of $\mathcal{L}$ with respect to the quasi-sure equivalence relation. We have the following results (see the Appendix).

THEOREM 2.1. *Each element in $\mathcal{L}$ can be identified with a quasi-continuous function on $\Omega$ (and so is defined quasi-surely). Moreover, $(L, c)$ is a Banach space.*

An element of $L$ is a class of equivalence, but, as usual, we do not heed the distinction between classes and their representatives. If $f \in \mathcal{L}$, there is a sequence $f_n \in C_b(\Omega)$ that converges to $f$ in $\mathcal{L}$. It is clear that for each $P \in \mathbf{P}$, $f_n$ form a Cauchy sequence in $L^2(\Omega, P)$ and, hence, converge to a function in $L^2(\Omega, P)$ equal to $f$ $P$-a.s.; this convergence is uniform in $P$. As a consequence, we get the following statement:

PROPOSITION 2.2. *Let $f \in \mathcal{L}$. Then $c(f) = \sup\{\|f\|_{L^2(\Omega,P)} : P \in \mathbf{P}\}$.*

REMARK. We consider here the capacity defined via $L^2$-norms (and not, e.g., via $L^1$-norms) only to make the stochastic calculus below easier. Nevertheless, if $f \in \mathcal{L}$, the quantity $\sup\{E_P|f| : P \in \mathbf{P}\}$ is well defined [and bounded by $c(f)$].

PROPOSITION 2.3. *Let $f \in C_b(\Omega)$, let $Q$ be a polynomial and let $s, t \in [0, T]$. Then the function $Q(B_t - B_s)f$ is in $\mathcal{L}$.*

PROOF. It is sufficient to consider the case where $Q(x) = x^k$, $k \geq 1$, and $s = 0$ (changes for arbitrary $s$ are obvious). Let $M$ be an upper bound of $|f|$. Using the Burkholder–Davis–Gundy inequalities and the hypothesis $H(\overline{\mu})$, we get that, for some constant $C_k > 0$,

$$\sup_{P \in \mathbf{P}} E_P |B_t|^{4k} \leq C_k \overline{\mu}_t^{2k} < \infty.$$

The function $S_n = (((-n) \vee B_t^k) \wedge n)f$ is in $C_b(\Omega)$. For any $P \in \mathbf{P}$, we have

$$\|B_t^k f - S_n\|_{L^2(\Omega,P)}^2 = E_P[\mathbf{1}_{\{|B_t^k| > n\}}(fB_t^k - S_n)^2]$$

$$\leq M^2 E_P[\mathbf{1}_{\{|B_t^k| > n\}} B_t^{2k}]$$

$$\leq M^2 P(|B_t^k| > n)^{1/2}(E_P B_t^{4k})^{1/2}$$

$$\leq M^2 \frac{E_P B_t^{4k}}{n^2} \leq C_k M^2 \frac{\overline{\mu}_t^{2k}}{n^2}.$$

The right-hand side here does not depend on $P$. Thus, $S_n$ is a Cauchy sequence in $\mathcal{L}$ that converges to $B_t^k f$ as $n \to \infty$ and the proof is complete. $\square$

### 2.2. Stochastic integrals in the presence of model uncertainty.
We define a stochastic integral suitable for modeling in our uncertain framework.

#### 2.2.1. Construction of the stochastic integral.
Let $\mathcal{H}_e$ be the set of "elementary" processes $h_s = \sum_{i=0}^N k_{t_i} \mathbb{1}_{]t_i, t_{i+1}]}(s)$, where $t_i$ is a deterministic subdivision of $[0, T]$ and $k_{t_i}$ are $\mathcal{F}_{t_i}$-measurable, bounded and continuous. We denote by $\mathcal{H}$ the completion of $\mathcal{H}_e$ with respect to the semi-norm

$$\|h\|_{\mathcal{H}} = c \left( \left( \int_0^T h_s^2 \, d\overline{\mu}_s \right)^{1/2} \right) = \sup_{P \in \mathbf{P}} E_P \left( \int_0^T h_s^2 \, d\overline{\mu}_s \right)^{1/2}$$

and denote by $H$ the quotient of $\mathcal{H}$ with respect to the linear space of processes $h$ such that $\|h\|_{\mathcal{H}} = 0$. It is clear that $H$ is a Banach space with respect to the resulting norm. Moreover, by the same type of arguments as in the case of real-valued functions, we get that each element $h$ in $\mathcal{H}$ (or in $H$) admits a predictable version and the mapping $\omega \mapsto h.(\omega)$ from $\Omega$ to $L^2([0, T], \overline{\mu})$ is quasi-continuous.

LEMMA 2.4. Let $h$ be a predictable process which admits a version in $L^2([0, T], \overline{\mu}; \mathcal{L})$. Then $h$ belongs to $\mathcal{H}$.

PROOF. Assume first that $h$ is in $\mathcal{H}_e$. Then we have

$$\|h\|_{\mathcal{H}}^2 \leq \sup_{P \in \mathbf{P}} E_P \int_0^T h_s^2 \, d\overline{\mu}_s \leq \int_0^T \sup_{P \in \mathbf{P}} E_P h_s^2 \, d\overline{\mu}_s = \int_0^T c^2(h_s) \, d\overline{\mu}_s.$$

Consider now $h \in L^2([0, T], \overline{\mu}; \mathcal{L})$ and assume that it is predictable. Then there exists a sequence in $\mathcal{H}_e$ which converges to $h$ in $L^2([0, T], \overline{\mu}; \mathcal{L})$. Now, the previous inequality ensures that it converges also in $\mathcal{H}$ to $h$ and we conclude. $\square$

As a corollary of this lemma and Proposition 2.3 we obtain the next lemma.

LEMMA 2.5. Let $Q$ be a polynomial. Then the process $Q(B)$ belongs to $\mathcal{H}$.

The following result follows from Lemma 2.4. We shall use it when $\varphi$ is an indicator function of a time interval.

LEMMA 2.6. *If $h \in L^2([0,T], \overline{\mu}; \mathcal{L})$ is predictable and $\varphi$ is a bounded deterministic function that is left-continuous with limits from the right, then the process $\varphi h$ belongs to $\mathcal{H}$.*

The next lemma will be useful when we will iterate the Itô formula.

LEMMA 2.7. *Let $h \in \mathcal{H}$. Then $\int_0^T h_s \, d\overline{\mu}(s)$ belongs to $\mathcal{L}$ and*

$$c\left(\int_0^T h_s \, d\overline{\mu}_s\right) \leq \overline{\mu}_T^{1/2} \|h\|_{\mathcal{H}}.$$

*As a consequence, the process $X = (X_t)$ with*

$$X_t = \int_0^t h_s \, d\overline{\mu}_s = \int_0^T h_s \mathbb{1}_{]0,t]}(s) \, d\overline{\mu}_s$$

*belongs to $\mathcal{H}$.*

The stochastic integral is defined in the following theorem.

THEOREM 2.8. *The linear mapping*

$$h = \sum_{i=0}^{N} k_{t_i} \mathbb{1}_{]t_i,t_{i+1}]} \mapsto I_T(h) = \int_0^T h_s \, dB_s = \sum_{i=0}^{N} k_{t_i}(B_{t_{i+1}} - B_{t_i}),$$

*considered as a function from $\mathcal{H}_e$ to $\mathcal{L}$, admits the bound*

(2) $$c(I_T(h)) \leq \|h\|_{\mathcal{H}}.$$

*It can be extended uniquely to a continuous linear mapping from $\mathcal{H}$ to $\mathcal{L}$, still denoted by $I_T(h) = \int_0^T h_s \, dB_s$, that satisfies* (2).

PROOF. We first assume that $h$ is an elementary process. Then, as a consequence of Proposition 2.3, $I_T(h)$ belongs to $\mathcal{L}$ and we have, for each $P \in \mathbf{P}$, that

$$\|I_T(h)\|_{L^2(\Omega,P)}^2 \leq E_P \int_0^T h_s^2 \, d\overline{\mu}_s.$$

Taking the supremum over all probabilities in $\mathbf{P}$, we get the desired inequality and conclude using a density argument. □

REMARK. Generalizing classical ideas, one can easily construct a stochastic integral with regular trajectories and other good properties such as Doob's inequality, see [9].

2.2.2. *Some properties of stochastic integrals.* Put $K = \{I_T(h) : h \in \mathcal{H}\}$. In the financial context the elements of this linear space are interpreted as the terminal values of portfolio processes. It is interesting to know whether $K$ is closed; we give a sufficient condition for this in the Appendix.

We turn now to estimates for powers of the canonical process.

PROPOSITION 2.9. *Let $s, t \in [0, T]$ and let $n$ be an integer. Then there exist $h \in \mathcal{H}$ and a positive constant $C$ depending on $n$ such that*

$$(B_t - B_s)^{2n} \leq \int_s^t h_u \, dB_u + C\overline{\mu}(]s, t])^n \qquad q.s.$$

PROOF. We give arguments for $s = 0$, because their extension to the general case is obvious. Let $P$ be in $\mathbf{P}$. By the Itô formula, we have

$$B_t^{2n} = 2n \int_0^t B_u^{2n-1} \, dB_u + n(2n-1) \int_0^t B_u^{2n-2} \, d\langle B \rangle_u^P$$

$$\leq 2n \int_0^t B_u^{2n-1} \, dB_u + n(2n-1) \int_0^t B_u^{2n-2} \, d\overline{\mu}_u, \qquad P\text{-a.s.}$$

For $n = 1$ the assertion follows from Lemmas 2.6 and A.7. For $n \geq 2$ we apply the Itô formula to $B_u^{2n-2}$ and the Fubini theorem to obtain that

$$\int_0^t B_u^{2n-2} \, d\overline{\mu}_u = (2n-2) \int_0^t \overline{\mu}(]u, t]) B_u^{2n-3} \, dB_u$$

$$+ (n-1)(2n-3) \int_0^t \overline{\mu}(]u, t]) B_u^{2n-4} \, d\langle B \rangle_u^P.$$

The integrator of the stochastic integral is an element of $\mathcal{H}$, while the ordinary integral admits the bound

$$\int_0^t \overline{\mu}(]u, t]) B_u^{2n-4} \, d\langle B \rangle_u^P \leq \overline{\mu}(]0, t]) \int_0^t B_u^{2n-4} \, d\overline{\mu}_u, \qquad P\text{-a.s.}$$

Continuing the reduction and applying Lemma A.7 at the end, we obtain the result. $\square$

The next lemma ensures that there is a "universal" version of the quadratic variation.

LEMMA 2.10. *Let $t \in [0, T]$. Then there exists an element in $\mathcal{L}$ that we denote by $\langle B \rangle_t$ such that $\langle B \rangle_t^P = \langle B \rangle_t$ $P$-a.s. for all $P \in \mathbf{P}$. Moreover,*

$$\langle B \rangle_t \leq \overline{\mu}_t \qquad q.s.$$

PROOF.   We just have to note that

$$\langle B \rangle_t = B_t^2 - 2 \int_0^t B_s \, dB_s.$$ □

A martingale that has a terminal value bounded from below by a constant $-a$ is bounded by this constant as a process almost surely. We have the following analog of this assertion, which follows directly from Lemma A.7.

LEMMA 2.11.   *Let $I_T(h) \in K$. Suppose that $I_T(h) > -a$ q.s. for some $a \geq 0$. Then $I_t(h) > -a$ q.s. for every $t \in [0, T]$.*

2.2.3. *Polar sets and stochastic integrals.*   It is natural to consider a stochastic integral under any probability $P' \in \mathbf{P}_m$ which does not charge polar sets (for the latter property, the notation $P' \ll c$ is natural). Observe that if $h \in \mathcal{H}$, then

$$\int_0^T h_u^2 \, d\overline{\mu}_u < \infty, \qquad P'\text{-a.s.}$$

So, by definition of a martingale probability, the process $\int_0^t h_u \, dB_u$ is a $P'$-local martingale and $\int_0^T h_s \, dB_s$ is well defined $P'$-a.s.

On the other hand, $\int_0^T h_s \, dB_s$ is defined q.s., so $P'$-a.s. it is defined as an element of $\mathcal{L}$, but by a density argument it is clear that these two definitions coincide. This has an interesting consequence: Any martingale probability which does not charge polar sets satisfies the same bracket assumption as the initial set of probabilities:

PROPOSITION 2.12.   *Assume $H(\underline{\mu}, \overline{\mu})$. Let $P' \in \mathbf{P}_m$ and $P' \ll c$. Then*

$$d\underline{\mu}_t \leq d\langle B \rangle_t^{P'} \leq d\overline{\mu}_t, \qquad P'\text{-a.s.}$$

PROOF.   Let $s \leq t \leq T$. Notice that the process $\int_s^{\cdot} (B_u - B_s) \, dB_u$, is well defined as a $P'$-local martingale on $[s, T]$. By the Itô formula,

$$(B_t - B_s)^2 = 2 \int_s^t (B_u - B_s) \, dB_u + \langle B \rangle_t^{P'} - \langle B \rangle_s^{P'}$$

$P'$-a.s. On the other hand,

$$(B_t - B_s)^2 \leq 2 \int_s^t (B_u - B_s) \, dB_u + \overline{\mu}([s, t])$$

quasi-surely, hence, $P'$-almost surely. This yields

$$\langle B \rangle_t^{P'} - \langle B \rangle_s^{P'} \leq \overline{\mu}([s, t]).$$

In the same way we have

$$(B_t - B_s)^2 \geq 2 \int_s^t (B_u - B_s) \, dB_u + \underline{\mu}([s,t])$$

quasi-surely and so we get the other inequality. □

2.3. *The superreplication problem.* In a financial context, an element $f$ of $\mathcal{L}$ can be interpreted as a contingent claim, that is, the cheapest riskless superreplication price in which we are interested. Put

$$\Lambda(f) = \inf\{a : \exists\, g \in K \text{ such that } a + g \geq f \text{ q.s.}\},$$

defining in this way a convex homogeneous function $\operatorname{dom}\Lambda = \{f \in \mathcal{L} : \Lambda(f) < \infty\}$. If $f \leq c = \text{const q.s.}$, then $\Lambda(f) \leq c$, in particular, with the usual identification, we have the inclusion $C_b(\Omega) \subset \operatorname{dom}\Lambda$.

2.3.1. *Properties of $\Lambda$.* The first result is a consequence of Proposition 2.9:

PROPOSITION 2.13. *For any integer $n \geq 1$, there exists a constant $C_{2n} > 0$ such that*

$$\Lambda((B_t - B_s)^{2n}) \leq C_{2n} \overline{\mu}([s,t])^n.$$

The second result deals with an approximation of the $B$. Let $t_i = it/n$ and let

$$S_t^n = \sum_{i=0}^{n-1} (B_{t_{i+1}} - B_{t_i})^2.$$

LEMMA 2.14. *Suppose that the distribution function $\overline{\mu}$ is Hölder continuous. Then*

$$\lim_{n \to \infty} \Lambda((S_t^n - \langle B \rangle_t)^2) = 0.$$

PROOF. By the Itô formula, we have that, for all $P \in \mathbf{P}$,

$$(S_t^n - \langle B \rangle_t)^2 = 8 \sum_{i=0}^{n-1} \int_{t_i}^{t_{i+1}} \left( \sum_{i=0}^{n-1} \int_{t_i \wedge s}^{t_{i+1} \wedge s} (B_u - B_{t_i}) \, dB_u \right) (B_s - B_{t_i}) \, dB_s$$

$$\qquad + 4 \sum_{i=0}^{n-1} \int_{t_i}^{t_{i+1}} (B_s - B_{t_i})^2 \, d\langle B \rangle_s^P$$

$$\leq 8 \sum_{i=0}^{n-1} \int_{t_i}^{t_{i+1}} \left( \sum_{i=0}^{n-1} \int_{t_i \wedge s}^{t_{i+1} \wedge s} (B_u - B_{t_i}) \, dB_u \right) (B_s - B_{t_i}) \, dB_s$$

$$\qquad + 4 \sum_{i=0}^{n-1} \int_{t_i}^{t_{i+1}} (B_s - B_{t_i})^2 \, d\overline{\mu}_s$$

$P$-a.s. Moreover, for all $u \leq s$,

$$(B_s - B_u)^2 = 2 \int_u^s (B_v - B_u) \, dB_v + \langle B \rangle_s^P - \langle B \rangle_u^P$$

$$\leq 2 \int_u^s (B_v - B_u) \, dB_v + \overline{\mu}([u, s]), \qquad P\text{-a.s.}$$

and, as in the proof of Proposition 2.9, we deduce that there exists $h \in \mathcal{H}$ such that

$$(S_t^n - \langle B \rangle_t)^2 \leq \int_0^t h_s \, dB_s + 4 \sum_{i=0}^{n-1} \int_{t_i}^{t_{i+1}} \overline{\mu}([s, t_i]) \, d\overline{\mu}_s$$

$$\leq \int_0^t h_s \, dB_s + 4C \left( \frac{t}{n} \right)^\alpha \overline{\mu}_t$$

quasi-surely, where $C$ and $\alpha$ are the constants from the Hölder condition (1). This yields

$$\Lambda((S_t^n - \langle B \rangle_t)^2) \leq 4C \left( \frac{t}{n} \right)^\alpha \overline{\mu}_t$$

and hence the result. $\square$

2.3.2. $\Lambda$ *and martingale measures.* We try to express $\Lambda(f)$ in terms of martingale measures which do not charge polar sets. The set of such measures we denote by $\mathbf{P}'$; $\mathcal{L}_+$ is the set of nonnegative (quasi-surely) functions in $\mathcal{L}$.

For any $P \in \mathbf{P}'$ and $f \in \mathcal{L}_+ \cap \mathrm{dom}\,\Lambda$, the function $f$ is defined $P$-a.s. Moreover, if $h \in \mathcal{H}$, the stochastic integral $\int_0^t h_u \, dB_u$ is a $P$-local martingale. Take $a \in \mathbb{R}$ and $g = \int_0^T h_s \, dB_s \in K$ such that $a + g \geq f$ q.s. Then for each $a + g \geq f$, $P$-a.s. is $P \in \mathbf{P}'$.

Observe now that thanks to Lemma 2.11, the stochastic integral $\int_0^t h_u \, dB_u$ is bounded from below. Since the process $\int_0^t h_s \, dB_s$ is a $P$-local martingale when $P \in \mathbf{P}'$, it follows from the Fatou lemma that in such a case $g \in L^1(P)$ and $E_P g \leq 0$. This yields that for any $P \in \mathbf{P}'$,

$$E_P[a + g] \geq E_P f$$

and, hence, $a \geq E_P f$. Summarizing:

LEMMA 2.15.   *The following holds:*

$$\mathcal{L}_+ \cap \mathrm{dom}\,\Lambda \subset \bigcap_{P \in \mathbf{P}'} L^1(P).$$

*Moreover, for* $f \in \mathcal{L}_+ \cap \mathrm{dom}\,\Lambda$,

(3)            $$\Lambda(f) \geq \sup\{E_P f : P \in \mathbf{P}'\} \geq \sup\{E_P f : P \in \mathbf{P}\}.$$

Our next goal is to establish the converse inequality for a suitable set of martingale probabilities.

**3. Main result.** The main result is established for contingent claims which belong to a large class $\Gamma \subset C_b(\Omega)$ (see Lemmas 5.4–5.6).

THEOREM 3.1. *Assume* $H(\underline{\mu}, \overline{\mu})$. *Then there exists a subset* $\mathbf{P}'' \subset \mathbf{P}_m$ *such that its elements also satisfy* $H(\underline{\mu}, \overline{\mu})$ *and for every* $f \in \Gamma$,

$$\Lambda(f) = \sup\{E_P f : P \in \mathbf{P}''\}.$$

Unfortunately, we do not know whether the set $\mathbf{P}''$ is the whole set of martingale measures that do not charge polar sets, because it is in the discrete-time case (cf. [10]). Nevertheless, this holds in the case when $\Lambda$ is finite (hence continuous) on $\mathcal{L}$.

3.1. *The particular case*: $\mathrm{dom}\,\Lambda = \mathcal{L}$.

PROPOSITION 3.2. *Assume* $\mathrm{dom}\,\Lambda = \mathcal{L}$. *Then for any* $f$ *in* $\mathcal{L}$,

$$(4) \qquad \Lambda(f) = \sup\{E_{P'} f : P' \in \mathbf{P}'\},$$

*where* $\mathbf{P}'$ *is the set of martingale measures that do not charge polar sets.*

PROOF. $\Lambda$ is well defined on $\mathcal{L}$ and, moreover, it is a sublinear map. As a consequence of the Hahn–Banach theorem, we have that

$$\Lambda(f) = \sup_{\lambda \in \mathcal{Q}} \lambda(f) \qquad \forall f \in \mathcal{L},$$

where $\mathcal{Q}$ is the set of linear mappings from $\mathcal{L}$ into $\mathbb{R}$ dominated by $\Lambda$.

LEMMA 3.3. *If* $\lambda \in \mathcal{Q}$, *then* $\lambda$ *is positive.*

PROOF. Let $f \in \mathcal{L}_+$. Then $\lambda(-f) \leq \Lambda(-f) \leq \Lambda(0) \leq 0$ and the result follows. □

LEMMA 3.4. *If* $\lambda \in \mathcal{Q}$, *then* $\lambda(1) = 1$ *and for all* $g \in K$, *we have* $\lambda(g) = 0$.

PROOF. The first assertion is obvious because $\lambda(1) \leq \Lambda(1) \leq 1$ and $\lambda(-1) \leq \Lambda(-1) \leq -1$. The second assertion follows because $\alpha\lambda(g) \leq \Lambda(\alpha g)$ for any real $\alpha$ and, as easily seen, $\Lambda(\alpha g) \leq 0$. □

It remains to check that the elements of $\mathcal{Q}$ are Borel measures that do not charge polar sets. This fact follows from a much more general result proved by Feyel and De La Pradelle (see [12], Proposition 11) which ensures that the dual space of $\mathcal{L}$ is a set of Borel measures that do not charge polar sets. So if $\lambda \in \mathcal{Q}$, there exists a measure $P'$ which belongs to $\mathbf{P}_m$ (as a consequence of Lemmas 3.3 and 3.4) such that

$$\forall f \in \mathcal{L} \qquad \lambda(f) = E_{P'}[f].$$

The proof of Proposition 3.2 is complete. $\square$

Unfortunately, the equality $\operatorname{dom}\Lambda = \mathcal{L}$ is not clear. In the rest of the paper, we try to characterize $\Lambda(f)$ in another way, namely, instead of the space $\mathcal{L}$, we work with continuous and bounded functions. To make use of standard representation theorems, we need to compactify $\Omega$. To simplify the proof, we consider first the case where all trajectories are bounded by a constant.

**4. Proof of the main result: Bounded case.** In this section we assume that $\Omega$ is a closed and bounded subset of $C([0,T];\mathbb{R})$, that is, the trajectories in $\Omega$ are bounded in absolute value by some $\gamma > 0$. Note that the integration theory above can be developed without any changes for such $\Omega$.

From now on, we assume that $\overline{\mu}$ is Hölder continuous.

4.1. *The Stone–Čech compactification of* $\Omega$. Let $\tilde{\Omega}$ denote a Stone–Čech compactification of the completely regular space $\Omega$ endowed with the supremum norm, that is, a pair $(\tilde{\Omega}, \phi)$, where $\tilde{\Omega}$ is compact and $\phi$ is a mapping $\phi : \Omega \to \tilde{\Omega}$ such that:

(i) $\phi$ is an homeomorphism onto $\phi(\Omega)$ equipped with the topology induced by $\tilde{\Omega}$.

(ii) $\phi(\Omega)$ is dense in $\tilde{\Omega}$.

(iii) For any bounded continuous function $f$ on $\Omega$ there exists a continuous function $\tilde{f}$ on $\tilde{\Omega}$ such that $f = \tilde{f} \circ \phi$.

Note that, conversely, for a continuous function $\tilde{f}$ on $\tilde{\Omega}$, the function $f = \tilde{f} \circ \phi$ is a bounded continuous function on $\Omega$. One can explicitly construct $\tilde{\Omega}$ and $\phi$ in the following way:

Let $C_\infty$ be the space of bounded complex functions on $\Omega$ considered as a commutative complex Banach algebra (cf. [3]). We shall choose the Stone–Čech compactification given by the character space of $C_\infty$ equipped with its weak topology (cf. [3]). In particular, $\phi$ is defined on $\Omega$ by

$$\phi(\omega)(x) = x(\omega), \qquad x \in C_\infty,$$

and for a function $f \in C_b(\Omega)$, its extensions $\tilde{f}$ is defined by

$$\tilde{f}(\tilde{\omega}) = \tilde{\omega}(f), \qquad \tilde{\omega} \in \tilde{\Omega}.$$

Since for every $t$, the function $B_t$ belongs to $C_\infty$, this entails the following information:

PROPOSITION 4.1. *Let $\tilde{B}_t$ denote the unique bounded continuous extension to $\tilde{\Omega}$ of the bounded continuous mapping $B_t$. Then*

$$\phi(\omega)(B_t) = B_t(\omega) = \tilde{B}_t(\phi(\omega)), \qquad \omega \in \Omega,$$

$$\tilde{\omega}(B_t) = \tilde{B}_t(\tilde{\omega}), \qquad \qquad \tilde{\omega} \in \tilde{\Omega}.$$

4.2. *Representation by probability measures on $\tilde{\Omega}$.* Let $\tilde{f}$ be a continuous (hence, bounded) function on $\tilde{\Omega}$. The function $f = \tilde{f} \circ \phi$ being bounded and continuous on $\Omega$ may be viewed as an element of $\mathcal{L}$. Consider the mapping

$$\tilde{\Lambda}: \tilde{f} \to \Lambda(f).$$

It is well defined on $C(\tilde{\Omega})$, the set of continuous functions on $\tilde{\Omega}$, and sublinear. By the Hahn–Banach theorem,

$$\tilde{\Lambda}(\tilde{f}) = \sup_{\lambda \in \mathcal{Q}} \lambda(\tilde{f}),$$

where $\mathcal{Q}$ is the set of linear mappings from $C(\tilde{\Omega})$ into $\mathbb{R}$ which are dominated by $\tilde{\Lambda}$. Since $\tilde{\Lambda}(0) = \Lambda(0) \le 0$, by the same proof as the one of Lemma 3.3, each $\lambda \in \mathcal{Q}$ is a positive linear form on $C(\tilde{\Omega})$. Therefore, since $\tilde{\Omega}$ is compact, it is a measure that we denote by $Q$. Now $Q(1) \le \tilde{\Lambda}(1) \le 1$ and $Q(-1) \le \tilde{\Lambda}(-1) \le -1$, so that $Q(1) = 1$ and $Q$ is a probability measure on $\tilde{\Omega}$. We want to show that each probability $Q \in \mathcal{Q}$ is a martingale measure. Let $\tilde{\mathcal{F}}_t = \sigma\{\tilde{B}_u, u \le t\}$.

LEMMA 4.2. *Let $Q \in \mathcal{Q}$. Then $\tilde{B}_t$ is a $Q$-martingale with respect to $\tilde{\mathcal{F}}_t$.*

PROOF. Let $s < t \le T$ and let $F: \mathbb{R}^d \mapsto \mathbb{R}$ be a bounded continuous function. Consider

$$f = F(B_{t_1}, \ldots, B_{t_d}),$$

where $t_i \le s$. We want to prove that $E_Q[(\tilde{B}_t - \tilde{B}_s)\tilde{f}] = 0$.

By standard arguments, $\Lambda((B_t - B_s)f) \le 0$ and $\Lambda(-(B_t - B_s)f) \le 0$, so that $E_Q \tilde{A} = 0$, where $A = (B_t - B_s)f$. However, in Lemma 4.6 below we check that $\tilde{A} = (\tilde{B}_t - \tilde{B}_s)F(\tilde{B}_{t_1}, \ldots, \tilde{B}_{t_d})$. This permits us to conclude. □

By applying the Cauchy–Schwarz inequality for each $Q \in \mathcal{Q}$ we get the following statement:

LEMMA 4.3. *Let $\tilde{f}$ and $\tilde{g}$ be in $C(\tilde{\Omega})$. Then*

$$\tilde{\Lambda}^2(\tilde{f}\tilde{g}) \le \tilde{\Lambda}(\tilde{f}^2)\tilde{\Lambda}(\tilde{g}^2).$$

Fix now $Q \in \mathcal{Q}$. We denote by $D$ the set of dyadic numbers of $[0, T]$.

PROPOSITION 4.4. *The function $t \in D \mapsto \tilde{B}_t(\tilde{\omega})$ is uniformly continuous for $Q$-almost all $\tilde{\omega} \in \tilde{\Omega}$, so it admits a continuous modification which we denote by $(\tilde{B}_t)_{t \in [0,T]}$. Moreover, if we assume $H(\underline{\mu}, \overline{\mu})$, then*

$$d\underline{\mu}_t \leq d\langle \tilde{B} \rangle_t^Q \leq d\overline{\mu}_t, \qquad Q\text{-}a.s.$$

PROOF. By Proposition 2.13, $\Lambda((B_t - B_s)^{2n}) \leq C_{2n}\overline{\mu}([s,t])^n$. Therefore,

$$E_Q(\tilde{B}_t - \tilde{B}_s)^{2n} \leq C_{2n}\overline{\mu}([s,t])^n \leq C^n|t-s|^{\alpha n},$$

where $C, \alpha > 0$ are the constants from (1). Now, taking $n$ large enough and following the classical proof (cf. [19], Chapter 2) of the Kolmogorov lemma, the first inequality yields that the function $t \in D \mapsto \tilde{B}_t(\tilde{\omega})$ is uniformly continuous and we can put $\tilde{B}_t = \lim_{u \in D, u \to t} \tilde{B}_u$.

For the second assertion, we define in the same way as in Lemma 2.14,

$$S^n = \sum_{i=0}^{n-1}(B_{s+i(t-s)/n} - B_{s+(i+1)(t-s)/n})^2.$$

As we shall see in the proof of Lemma 4.6,

$$\tilde{S}^n = \sum_{i=0}^{n-1}(\tilde{B}_{s+i(t-s)/n} - \tilde{B}_{s+(i+1)(t-s)/n})^2, \qquad Q\text{-a.e.}$$

Let $f$ be a nonnegative function in $C_b(\Omega)$. We have

$$\lim_{n \to \infty} E_Q[\tilde{f}\{\tilde{S}^n - \overline{\mu}([s,t])\}] \leq \lim_{n \to \infty} \Lambda(f\{S^n - \overline{\mu}([s,t])\})$$
$$\leq \lim_{n \to \infty} \Lambda(f\{S^n - (\langle B \rangle_t - \langle B \rangle_s)\})$$
$$+ \Lambda(f\{\langle B \rangle_t - \langle B \rangle_s - \overline{\mu}([s,t])\}).$$

Whereas $\langle B \rangle_t - \langle B \rangle_s - \overline{\mu}([s,t]) \leq 0$ q.s. and $f$ is nonnegative,

$$\Lambda(f\{\langle B \rangle_t - \langle B \rangle_s - \overline{\mu}([s,t])\}) \leq 0.$$

Thanks to Lemma 4.3,

$$\Lambda^2(f\{S^n - (\langle B \rangle_t - \langle B \rangle_s)\}) \leq \Lambda(f^2)\Lambda((S^n - (\langle B \rangle_t - \langle B \rangle_s))^2)$$

and by adapting the proof of Lemma 2.14, we obtain that

$$\lim_{n \to \infty} \Lambda((S^n - (\langle B \rangle_t - \langle B \rangle_s))^2) = 0.$$

So, passing to the limit, we get

$$E_Q[f\{\langle \tilde{B} \rangle_t - \langle \tilde{B} \rangle_s - \overline{\mu}([s,t])\}] \leq 0.$$

Using a density argument, we conclude that

$$\langle \tilde{\tilde{B}} \rangle_t - \langle \tilde{\tilde{B}} \rangle_s - \overline{\mu}([s,t]) \leq 0, \qquad Q\text{-a.s.}$$

This implies that the quadratic variation of $B$ under $Q$ is absolutely continuous with respect to $\overline{\mu}$ and is dominated by $\overline{\mu}$.

For the other inequality, a similar proof works if we consider $f \in C_b(\Omega)$ to be nonpositive and use the fact that

$$\langle B \rangle_t - \langle B \rangle_s - \underline{\mu}([s,t]) \geq 0 \qquad \text{q.s.} \qquad \square$$

We denote by $Q^*$ the law of $\tilde{\tilde{B}}$ under $Q$. It is clear that $Q^*$ is a martingale measure on $\Omega$. An immediate consequence of the preceding theorem is the next statement.

COROLLARY 4.5. *The process $(B_t)$ is a $Q^*$-martingale which satisfies the same hypothesis $H(\overline{\mu})$ or $H(\underline{\mu}, \overline{\mu})$ with which we begin.*

What is the relationship between $Q$ and $Q^*$? We claim that for a large class of functions $f$ in $C_b(\Omega)$,

$$(5) \qquad E_Q \tilde{f} = E_{Q^*} f.$$

Let $\Gamma$ be the set of bounded continuous functions on $\Omega$ that satisfy (5) for each $Q \in \mathcal{Q}$. The next lemmas show that standard options (with deterministic maturity $T$) belong to $\Gamma$.

LEMMA 4.6. *If $f$ is a cylindrical continuous function, it belongs to $\Gamma$.*

PROOF. We want to show that $f = F(B_{t_1}, \ldots, B_{t_d})$ belongs to $\Gamma$ when $F$ is a bounded continuous function. We have that

$$\tilde{f}(\tilde{\omega}) = \tilde{\omega}(f), \qquad \tilde{\omega} \in \tilde{\Omega}.$$

Let $\omega \in \Omega$. Then

$$\begin{aligned}
\tilde{f}(\phi(\omega)) &= \phi(\omega)(f) \\
&= f(w) \\
&= F(B_{t_1}(\omega), \ldots, B_{t_d}(\omega)) \\
&= F(\tilde{B}_{t_1}(\phi(\omega)), \ldots, \tilde{B}_{t_d}(\phi(\omega))).
\end{aligned}$$

So we have proved that

$$\tilde{f}(\tilde{\omega}) = F(\tilde{B}_{t_1}(\tilde{\omega}), \ldots, \tilde{B}_{t_d}(\tilde{\omega})), \qquad \tilde{\omega} \in \phi(\Omega).$$

Whereas $\phi(\Omega)$ is dense in $\tilde{\Omega}$ and each side of the previous equality is continuous,

$$\tilde{f}(\tilde{\omega}) = F(\tilde{B}_{t_1}(\tilde{\omega}), \ldots, \tilde{B}_{t_d}(\tilde{\omega})), \qquad \tilde{\omega} \in \tilde{\Omega}.$$

It is now easy to conclude.  $\square$

Another example provides functions that depend on the supremum of trajectory. Put

$$S = \sup_{t \in [0,T]} B_t = \sup_{t \in D} B_t.$$

PROPOSITION 4.7.  *Let* $G \colon \mathbb{R} \longrightarrow \mathbb{R}$ *be a continuous function. Then* $f = G(S)$ *is in* $\Gamma$.

PROOF.  We split the arguments into several steps:

*Step* 1. Let $\beta > 0$ and let $k$ be an integer such that $k\alpha - \beta > -1/2$, where $\alpha$ is the constant in (1). Consider now the everywhere defined mapping $Y \colon \Omega \to \mathbb{R} \cup \{\infty\}$:

$$Y = \int_0^T \int_0^T \frac{|B_t - B_s|^{2k}}{|t-s|^{\beta}} \, ds \, dt.$$

We claim that $Y$ belongs to $\mathcal{L}$. Indeed, using Lemma 2.9, we have for $P \in \mathbf{P}$

$$E_P Y^2 \leq T^2 \int_0^T \int_0^T \frac{E_P |B_t - B_s|^{4k}}{|t-s|^{2\beta}} \, ds \, dt$$

$$\leq \int_0^T T^2 C_{4k} \int_0^T \frac{|t-s|^{2\alpha k}}{|t-s|^{2\beta}} \, ds \, dt$$

$$= C',$$

where $C'$ is a finite constant which only depends on $\alpha$, $\beta$ and $k$, and not on $P$. So we have

$$\sup_{P \in \mathbf{P}} E_P(Y^2) < \infty.$$

Put

$$Y_n = \int_0^T \int_0^T \frac{|B_t - B_s|^{2k}}{|t-s|^{\beta} + 1/n} \, ds \, dt.$$

It is clear that $Y_n$ belongs to $\mathcal{L}$ (because it is bounded and continuous) and estimates similar to the preceding one yield that $Y_n$ is a Cauchy sequence in $\mathcal{L}$ and hence converges to $Y$. As a consequence, $Y(w)$ is finite q.s. and the so-called Garsia–Rademich–Rumsey inequality (see [15]) ensures that if one takes $\beta > 2$ and $k$ large enough so that $\alpha k - \beta > -1/2$, then there exists a constant $C''$ which only depends on $\beta$ and $k$ such that quasi-surely

$$(6) \qquad |B_u(w) - B_v(w)|^{2k} \leq C'' Y(w)|u - v|^{\beta - 2} \qquad \forall u, v \in [0, T].$$

*Step* 2. We now prove that $Y$ belongs to dom $\Lambda$. As in the proof of Proposition 2.9, we have the bound

$$(B_t - B_s)^{2k} \leq \int_s^t h_u^{s,t} \, dB_u + C_{2k}\overline{\mu}([s, t])^k,$$

where $h^{s,t}$ is such that, for $u \in \, ]s, t]$,

$$h_u^{s,t} = 2k(B_u - B_s)^{2k-1} + \sum_{i=1}^{k-1} \gamma_i (\overline{\mu}([s, t]))^{i-1} \overline{\mu}([u, t])(B_u - B_s)^{2k-1-2i},$$

where $\gamma_i$ are some constants that depend on $k$. If $u \notin \, ]s, t]$, we set $h_u^{s,t} = 0$. We have

$$\sup_{P \in \mathbf{P}} \int_0^T E_P \left( \int_0^T \int_0^t \frac{h_u^{s,t}}{|t - s|^\beta} \, ds \, dt \right)^2 d\overline{\mu}_u \leq M \int_0^T \int_0^t \frac{\overline{\mu}([s, t])^{2k}}{|t - s|^{2\beta}} \, ds \, dt < \infty,$$

where $M$ is some constant. From this, by an approximation argument, it is clear that the process

$$h : u \to \int_0^T \int_0^t \frac{h_u^{s,t}}{|t - s|^\beta} \, ds \, dt$$

belongs to $L^2([0, T], \overline{\mu}; \mathcal{L})$ and hence to $\mathcal{H}$. Moreover, we have that

$$Y = 2 \int_0^T \int_0^t \frac{|B_t - B_s|^{2k}}{|t - s|^\beta} \, ds \, dt \leq 2 \int_0^T h_s \, dB_s + 2C_{2k}\zeta,$$

where

$$\zeta = \int_0^T \int_0^t \frac{\overline{\mu}([s, t])^k}{|t - s|^\beta} \, ds \, dt \leq C^k \int_0^T \int_0^t \frac{|t - s|^{\alpha k}}{|t - s|^\beta} \, ds \, dt < \infty,$$

because $\alpha k - \beta > -1/2$, so $Y$ belongs to dom $\Lambda$.

*Step* 3. We denote by $D_n$ the set of dyadic numbers of order $n$ in $[0, T]$,

$$D_n = \left\{ 0, \frac{T}{2^n}, \dots, \frac{T-1}{2^n}, T \right\},$$

so that $D = \bigcup_n D_n$. We put

$$S_n = \sup_{t \in D_n} B_t.$$

As a consequence of Lemma 4.6,

$$\tilde{S}_n(\tilde{\omega}) = \sup_{t \in D_n} \tilde{B}_t(\tilde{\omega}), \qquad \tilde{\omega} \in \tilde{\Omega}.$$

On the other hand, thanks to (6),

$$|S - S_n|^{2k} \leq C''Y \left(\frac{T}{2^n}\right)^{\beta-2} \qquad \text{q.s.}$$

Whereas $Y \in \operatorname{dom}\Lambda$, we get that

$$\lim_{n\to\infty} \Lambda((S - S_n)^{2k}) = 0$$

and so

$$\lim_{n\to\infty} E_Q[(\tilde{S} - \tilde{S}_n)^{2k}] = 0.$$

From this, we deduce that $\tilde{S}_n$ increases $Q$-a.s. to $\tilde{S}$ and

$$E_Q G(\tilde{S}) = \lim_{n\to\infty} E_Q G(\tilde{S}_n) = \lim_{n\to\infty} E_{Q^*} G(S_n) = E_{Q^*} G(S)$$

by virtue of the dominated convergence theorem. □

In the same way we get the following proposition.

PROPOSITION 4.8. *Let* $F, G: \mathbb{R} \longrightarrow \mathbb{R}$ *be two continuous functions. Then* $f = G(\int_0^T F(B_s)\,ds)$ *belongs to* $\Gamma$.

PROOF. First, we note that $f \in C_b(\Omega)$. Then we consider $\tilde{\omega} \in \Omega$ such that the function $t \in D \to \tilde{B}_t(\tilde{\omega})$ is uniformly continuous and, by construction, $t \in [0, T] \to \tilde{\tilde{B}}_t(\tilde{\omega})$ is continuous. Next, we define the sets

$$O_n^1 = \bigcap_{s \in D_n} \{w' \in \tilde{\Omega}, |w'(B_s) - \tilde{\omega}(B_s)| < 1/n\},$$

$$O_n^2 = \bigcap_{s \in D_n \setminus \{T\}} \left\{w' \in \tilde{\Omega}, w'\left(\sup_{u \in [0, T/2^n]} |B_{s+u} - B_s|\right) < 1/n\right\}$$

and, finally,

$$O_n = O_n^1 \cap O_n^2 \cap \{w' \in \tilde{\Omega}, |w'(f) - \tilde{\omega}(f)| < 1/n\}.$$

By definition of the weak topology, $O_n$ is an open set in $\tilde{\Omega}$. Because $\phi(\Omega)$ is dense, for each $n$ there exists $w_n \in \Omega$ such that $\Phi(w_n) \in O_n$. Now it is clear that for all $s \in D$,

$$\tilde{B}_s(\tilde{\omega}) = \lim_{n\to\infty} B_s(w_n).$$

For all $t \in [0, T]$ and all $n$, there is $t_n \in D_n$ such that $|t - t_n| < \frac{T}{2^n}$. Then, because $w_n \in O_n^2$,

$$|B_t(w_n) - B_{t_n}(w_n)| < 1/n,$$

which yields easily that

$$\tilde{\tilde{B}}_t(\tilde{\omega}) = \lim_{n \to \infty} B_t(w_n).$$

Moreover, by definition of $O_n$,

$$\tilde{f}(\tilde{\omega}) = \lim_{n \to \infty} f(w_n).$$

As a consequence of the dominated convergence theorem,

$$\tilde{f}(\tilde{\omega}) = \lim_{n \to \infty} f(w_n) = \lim_{n \to \infty} G\left(\int_0^T F(B_t(w_n))\, dt\right) = G\left(\int_0^T F(\tilde{\tilde{B}}_t(\tilde{\omega}))\, dt\right).$$

From this, we deduce that

$$\tilde{f} = G\left(\int_0^T F(\tilde{\tilde{B}}_t)\, dt\right)$$

$Q$-a.s. and the result follows. $\quad\square$

PROOF OF THEOREM 3.1 (Bounded case). Clearly, the theorem holds with $\mathbf{P}'' = \{Q^* : Q \in \mathcal{Q}\}$. $\quad\square$

## 5. Proof of the main result: Unbounded case.
We now apply the same method to the whole space $\Omega = C([0, T], \mathbb{R})$. We consider the same canonical Stone–Čech compactification as in Section 4.1 (also with the same notation). As in Section 4.2, $\tilde{\Lambda}$ is well defined with

$$\tilde{\Lambda}(\tilde{f}) = \sup_{\lambda \in \mathcal{Q}} \lambda(\tilde{f}), \qquad \tilde{f} \in C(\tilde{\Omega}),$$

where $\lambda \in \mathcal{Q}$ may be represented by a probability $Q$ on $\tilde{\Omega}$.

5.1. *Study of the process $\tilde{B}$.* Some care is needed since now $B_t$, for a given $t$, is no longer a bounded continuous function on $\Omega$. We define $\tilde{B}_t$ via a limiting procedure. To this end we need a lemma:

LEMMA 5.1. *Let $f \geq 0$ be a continuous function on $\Omega$. Then $f$ can be extended to a function $\tilde{f}$ on $\tilde{\Omega}$ with values in $[0, \infty]$.*

PROOF. The function $f \wedge n$ is bounded and continuous on $\Omega$ and, hence, admits a bounded continuous extension to $\tilde{\Omega}$ that we denote by $\tilde{f}_n$. For every $\tilde{\omega} \in \tilde{\Omega}$, $\tilde{f}_n(\tilde{\omega})$ is nondecreasing [since it is on the dense subset $\Phi(\Omega)$]. We set

$$\tilde{f}(\tilde{\omega}) = \lim_{n \to \infty} \tilde{f}_n(\tilde{\omega}), \qquad \tilde{\omega} \in \tilde{\Omega}. \qquad\square$$

If $f$ is a continuous function in $\Omega$, we set $\tilde{f^+} = \lim_{n \to \infty} \widetilde{f^+ \wedge n}$ and $\tilde{f^-} = \lim_{n \to \infty} \widetilde{f^- \wedge n}$.

*Definition and notation.* Let $f$ be a continuous function on $\Omega$. It may be extended to an everywhere defined function $\tilde{f}$ on $\tilde{\Omega}$ with values in $[-\infty, \infty]$ given, for $\tilde{\omega} \in \tilde{\Omega}$, by

$$\tilde{f}(\tilde{\omega}) = \begin{cases} \tilde{f^+}(\tilde{\omega}) - \tilde{f^-}(\tilde{\omega}), & \text{if } \tilde{f^+}(\tilde{\omega}) < \infty \text{ and } \tilde{f^-}(\tilde{\omega}) < \infty, \\ \infty, & \text{if } \tilde{f^+}(\tilde{\omega}) = \infty \text{ and } \tilde{f^-}(\tilde{\omega}) < \infty, \\ -\infty, & \text{if } \tilde{f^+}(\tilde{\omega}) < \infty \text{ and } \tilde{f^-}(\tilde{\omega}) = \infty, \\ \infty, & \text{if } \tilde{f^+}(\tilde{\omega}) = \infty \text{ and } \tilde{f^-}(\tilde{\omega}) = \infty. \end{cases}$$

As a consequence of this definition, the process $\tilde{B}_t$ is well defined.

Let $Q$ be a probability associated to an element $\lambda \in \mathcal{Q}$. By the Fatou lemma,

$$E_Q[(\tilde{B}_t^+)^2] \leq \liminf_{n \to \infty} E_Q[\widetilde{(B_t^+ \wedge n)^2}] \leq \liminf_{n \to \infty} \Lambda((B_t^+ \wedge n)^2) \leq \Lambda(B_t^2) \leq \overline{\mu}_t.$$

By the same reasoning for the negative part $B_t^-$, we obtain that $\tilde{B}_t$ belongs to $L^2(Q)$. In particular, it is $Q$-a.s. finite and, therefore,

$$\tilde{B}_t = \lim_{n \to \infty} (\widetilde{B_t^+ \wedge n} - \widetilde{B_t^- \wedge n})$$

is in $L^2(Q)$.

In the same way we can show that

$$E_Q[(\tilde{B}_t - \tilde{B}_s)^2] \leq \Lambda((B_t - B_s)^2)$$

and also that, for any $n$,

$$(7) \qquad E_Q[(\tilde{B}_t - \tilde{B}_s)^{2n}] \leq \Lambda((B_t - B_s)^{2n}) \leq C_{2n}\overline{\mu}([s,t])^n.$$

Now, we show that the process $(\tilde{B}_t)$ is a $Q$-martingale.

5.2. *The martingale property.* We cannot argue as in the bounded case because $\tilde{B}$ is not bounded, and $\tilde{\Lambda}(\tilde{A}(\tilde{B}_t - \tilde{B}_s))$ for a bounded continuous and $F_s$-measurable function $A$ is not well defined. We shall use instead an approximation.

Let us introduce the sequence $f_n$ of real-valued functions

$$f_n(x) = \begin{cases} x, & \text{if } x \in [-n, n], \\ n + n \arctan\left(\dfrac{x - n}{n}\right), & \text{if } x > n, \\ -n + n \arctan\left(\dfrac{x + n}{n}\right), & \text{if } x < n. \end{cases}$$

One can easily verify that $f_n$ is a $C^2$ function with bounded first and second derivatives and that $\lim_n f_n(x) = x$, $|f'(x)| \leq 1$ and $|f_n''(x)| \leq 1/n$. Recall that we assume that $\overline{\mu}$ is Hölder continuous.

PROPOSITION 5.2. *The process $\tilde{B}_t$ is a Q-martingale. Moreover, it admits a continuous modification that we denote by $\tilde{\tilde{B}}_t$.*

PROOF. Take a nonnegative cylindrical and continuous function $A$ which is $F_s$-measurable and bounded by $M > 0$. By the Itô formula, under any $P \in \mathbf{P}$, we have

$$Af_n(B_t - B_s) = \int_s^t Af'_n(B_u - B_s)\, dB_u + \frac{1}{2}\int_s^t Af''_n(B_u - B_s)\, d\langle B\rangle_u^P$$

so that

$$\Lambda(Af_n(B_t - B_s)) \le \frac{M}{2n}\overline{\mu}([s,t]).$$

If we set $G = Af_n(B_t - B_s)$, as in the proof of Lemma 4.6, we get

$$\tilde{G} = \tilde{A}f_n(\tilde{B}_t - \tilde{B}_s).$$

So we have that

$$E_Q[\tilde{A}f_n(\tilde{B}_t - \tilde{B}_s)] \le \frac{M}{2n}\overline{\mu}([s,t]).$$

In the same way, starting from $-A \cdot f_n(B_t - B_s)$ we get

$$-E_Q[\tilde{A}f_n(\tilde{B}_t - \tilde{B}_s)] \le \frac{M}{2n}\overline{\mu}([s,t]),$$

whence

$$\lim_{n\to\infty} E_Q[\tilde{A}f_n(\tilde{B}_t - \tilde{B}_s)] = 0.$$

It remains to show that the sequence $(\tilde{A}f_n(\tilde{B}_t - \tilde{B}_s))$ tends in $L^1(Q)$ to $\tilde{A}(\tilde{B}_t - \tilde{B}_s)$. However, this is clear, since

$$E_Q|\tilde{A}f_n(\tilde{B}_t - \tilde{B}_s) - \tilde{A}(\tilde{B}_t - \tilde{B}_s)| \le E_Q[|A\|\tilde{B}_t - \tilde{B}_s|\mathbb{1}_{|\tilde{B}_t - \tilde{B}_s| > n}]$$
$$\le M\frac{E_Q|\tilde{B}_t - \tilde{B}_s|^2}{n}$$
$$\le M\frac{C_2\overline{\mu}([s,t])}{n}.$$

Thus

$$E_Q[\tilde{A}(\tilde{B}_t - \tilde{B}_s)] = 0.$$

The last part of the assertion follows from (7) as in Proposition 4.4. □

Let us turn now to the quadratic variation of $\tilde{\tilde{B}}$.

LEMMA 5.3. *Assume $H(\underline{\mu}, \overline{\mu})$ and let $Q \in \mathcal{Q}$. Then*

$$d\underline{\mu}_t \leq d\langle \tilde{\tilde{B}} \rangle_t^Q \leq d\overline{\mu}_t, \qquad Q\text{-}a.s.$$

PROOF. We adopt the same notation as in the proof of Proposition 4.4, so we fix $s < t$. The same argument works except that the functions $S^n$ are no longer bounded. Nevertheless they are in $\operatorname{dom} \Lambda$ (defined in Section 2.3) and still the argument holds because we have

$$E_Q \tilde{S}^n \leq \Lambda(S^n).$$

To prove this, consider the sequence of real functions $f_k$ introduced at the beginning of this subsection. Clearly, $f_k(\tilde{S}^n) = \widetilde{f_k(S^n)}$ and $\lim_k E_Q(f_k(\tilde{S}^n)) = E_Q(\tilde{S}^n)$, so

$$E_Q \tilde{S}^n \leq \lim_{k \to \infty} \Lambda(f_k(S^n)) \leq \Lambda(S^n). \qquad \square$$

From now on, we can use the same arguments as in Section 4.2 to conclude. So, we still define $Q^*$ to be the law [on $(\Omega, \mathcal{B})$] of $(\tilde{\tilde{B}}_t)_{t \in [0,T]}$ and define by $\Gamma$ the set of bounded continuous functions $f$ such that

$$E_Q \tilde{f} = E_{Q^*} f.$$

The same proofs as in Section 4.2, with truncation arguments, give the following lemmas.

LEMMA 5.4. *If $f$ is a continuous and bounded cylindrical function, then $f \in \Gamma$.*

LEMMA 5.5. *Let $F : \mathbb{R} \longrightarrow \mathbb{R}$ be a continuous function and let $G : \mathbb{R} \longrightarrow \mathbb{R}$ be a bounded continuous function. Then $f = G(\int_0^T F(B_s)\,ds)$ is in $\Gamma$.*

Finally, we still denote by $S$ the function $S = \sup_{t \in [0,T]} B_t = \sup_{t \in D} B_t$. Then we can state the next lemma.

LEMMA 5.6. *Let $G : \mathbb{R} \longrightarrow \mathbb{R}$ be a bounded continuous function. Then $f = G(S)$ is in $\Gamma$.*

PROOF OF THEOREM 3.1 (Unbounded case). Here again we just have to put $\mathbf{P}'' = \{Q^* : Q \in \mathcal{Q}\}$. $\square$

**6. Application to a generalized UVM model.** Now the following theorem is a straightforward consequence of Theorem 3.1 and Lemma 2.15:

THEOREM 6.1. *Let $\underline{\mu}$ and $\overline{\mu}$ be two deterministic measures on $[0,T]$ such that $d\underline{\mu} \leq d\overline{\mu}$ and $\overline{\mu}$ is Hölder continuous. Let $\mathbf{P}$ be the set of all the martingale measures which satisfy the assumption $H(\underline{\mu}, \overline{\mu})$. Let $f$ be a bounded continuous function in $\Gamma$. Then*

$$\Lambda(f) = \sup\{E_P f : P \in \mathbf{P}\}.$$

Note that the UVM model corresponds to the case of the Lebesgue measure. As mentioned earlier, this result is new even in this case because it encompasses quite general path-dependent European options.

**7. Conclusion.** In this paper we set a framework for dealing with model uncertainty in the pricing of contingent claims. We provide a refined version of the stochastic integral which is defined quasi-surely with respect to a family of martingale laws on the canonical space which may not be dominated in the statistical sense. We study then the problem of the cheapest super-replication strategy. In the case when the bracket of the canonical process under the laws of the family is controlled, we give a partial characterization of the cheapest superhedging price by using a compactification method. In the case of the UVM model, this characterization is complete and it works for a large class of European path-dependent claims. The characterization in the general setting remains an open question.

## APPENDIX

Let $\mathbf{P}$ be a set of probability measures on the path space $\Omega$. For $f \in C_b(\Omega)$ we define

$$c(f) = \sup\{\|f\|_{L^2(\Omega,P)} : P \in \mathbf{P}\}.$$

Clearly, $c(f)$ is a semi-norm such that $c(1) = 1$; $c(f) = c(|f|)$, hence, $c(f) \leq 1$ if $|f| \leq 1$.

In a classical way [5, 8], we consider the Lebesgue extension of $c$:

- for lower semicontinuous $f \geq 0$,

$$c(f) = \sup\{c(\varphi) : \varphi \in C_b(\Omega), 0 \leq \varphi \leq f\};$$

- for arbitrary $g : \Omega \to \bar{\mathbb{R}}$,

$$c(g) = \inf\{c(f) : f \text{ is lower semicontinuous}, f \geq |g|\}.$$

For $A \subset \Omega$ we put $c(A) = c(\mathbf{1}_A)$.

The theory goes well under the following regularity property of $c$:

HYPOTHESIS (R).  $c(f_n) \downarrow 0$ for every sequence $f_n \in C_b(\Omega)$ such that $f_n \downarrow 0$.

THEOREM A.1.  *If (R) holds, the set function $c(A)$ is a (regular) Choquet capacity, that is, the following statements hold:*

1. $\forall A \in \mathcal{B},\ 0 \leq c(A) \leq 1$.
2. *If $A \subset B$, then $c(A) \leq c(B)$.*
3. *If $A_n$ is a sequence of sets in $\mathcal{B}$, then $c(\bigcup_n A_n) \leq \sum_n c(A_n)$.*
4. *If $A_n$ is an increasing sequence of sets in $\mathcal{B}$, then $c(\bigcup_n A_n) = \lim_n c(A_n)$.*

For any $f_n \in C_b(\Omega)$ the function $F_n: P \mapsto \|f_n\|_{L^2(\Omega, P)}$ is continuous with respect to the weak convergence of probability measures; if $f_n \downarrow 0$, then $F_n(P) \downarrow 0$ for any $P$. By the Dini lemma this convergence is uniform on weakly compact subsets and we get the next lemma:

LEMMA A.2.  *If $\mathbf{P}$ is relatively weakly compact (i.e., tight), then $c$ is regular.*

Let $c$ denote the capacity associated to $\mathbf{P}$ as defined in Section 2.1.1. The Rebolledo criterion (see Theorem VI.4.13 in [16]) says that $\mathbf{P}$ (the set of laws of continuous martingales) is tight if and only if the set of laws of $\langle B \rangle^P$, $P \in \mathbf{P}$, is tight. The latter property holds when $\mathbf{P}$ satisfies $H(\overline{\mu})$. We summarize this next:

LEMMA A.3.  *Under $H(\overline{\mu})$, Hypothesis (R) is satisfied, that is, $c$ is regular.*

We now study $\mathcal{L}$. Clearly, $\mathcal{L}$ contains $C_b(\Omega)$. In the converse direction, it is interesting to know that the analog of the Lusin theorem holds in our setting. The proof relies on the following simple fact.

LEMMA A.4.  *Let $f \in C_b(\Omega)$. Then for each $\alpha > 0$,*

$$c(\{|f| > \alpha\}) \leq \frac{c(f)}{\alpha}.$$

PROOF.  Whereas $f$ is continuous, $\mathbb{1}_{\{|f| > \alpha\}}$ is lower semicontinuous. Take arbitrary $\varphi \in C_b(\Omega)$ such that $0 \leq \varphi \leq \mathbb{1}_{\{|f| > \alpha\}}$. By the Markov inequality for any probability $P \in \mathbf{P}$, we have

$$(8) \qquad \|\varphi\|^2_{L^2(\Omega, P)} \leq P(|f| > \alpha) \leq \frac{\|f\|^2_{L^2(\Omega, P)}}{\alpha^2}.$$

Taking the supremum over $P$ and $\varphi$, we get the result.  □

Due to the $\sigma$-subadditivity of $c$, we verify that it satisfies a Borel–Cantelli lemma, which yields the next result (the proof can be found in [4]):

PROPOSITION A.5. *Let $f_n$ be a c-Cauchy sequence in $C_b(\Omega)$. Then, for each $\varepsilon > 0$, there exists an open set $O$ with $c(O) < \varepsilon$ such that $f_n$ converges uniformly on $O^c$.*

The following lemma, which is a consequence of this proposition, provides a bridge between the space $\mathcal{L}$ and the $L^2(P)$ spaces for $P \in \mathbf{P}$:

LEMMA A.6. *Let $f, g \in \mathcal{L}$ be such that $f = g$ $P$-a.s. for every $P \in \mathbf{P}$. Then $f = g$ in $\mathcal{L}$.*

LEMMA A.7. *Let $f, g \in \mathcal{L}$ be such that $f \leq g$ $P$-a.s. for every $P \in \mathbf{P}$. Then $f \leq g$ quasi-everywhere.*

PROOF. Let $h \in \mathcal{L}$ and let $h_n$ be a sequence in $C_b(\Omega)$ convergent to $h$ in $\mathcal{L}$. Using the inequality $|a^+ - b^+| \leq |a - b|$, we have, for any $P \in \mathbf{P}$,

$$\|h_n^+ - h_m^+\|_{L^2(\Omega, P)} \leq \|h_n - h_m\|_{L^2(\Omega, P)} \leq c(h_n - h_m).$$

Taking the supremum over $P \in \mathbf{P}$, we conclude that $h_n^+$ is a Cauchy sequence in $\mathcal{L}$ which clearly converges to $h^+$. Thus, $h^+$ belongs to $\mathcal{L}$. To conclude, we apply the previous lemma to the functions $(f - g)^+$ and $0$. □

We end this Appendix with a sufficient condition that ensures that the subspace of stochastic integrals $K = \{I_T(h) : h \in \mathcal{H}\}$ is closed in $\mathcal{L}$:

THEOREM A.8. *Assume that $H(a\overline{\mu}, \overline{\mu})$ holds with some $a \in \,]0, 1[$. Then $K$ is closed.*

PROOF. Let $f_n = I_T(h^n)$, $h^n \in \mathcal{H}$, form a sequence which converges to $f$ in $\mathcal{L}$. For any $P \in \mathbf{P}$, we have the inequality

$$E_P\left(\int_0^T (h_s^n - h_s^m)^2 \, d\overline{\mu}_s\right)^{1/2} \leq \frac{1}{a^{1/2}} \|f_n - f_m\|_{L^2(\Omega, P)}.$$

Taking the supremum over $P$, we get that

$$c\left(\left(\int_0^T (h_s^n - h_s^m)^2 \, d\overline{\mu}_s\right)^{1/2}\right) \leq \frac{1}{a^{1/2}} c(f_n - f_m),$$

that is, $h^n$ is a Cauchy sequence in $\mathcal{H}$, hence, it converges to a limit $h \in \mathcal{H}$. It is easy to verify that $f = I_T(h)$. □

**Acknowledgment.** The authors are very grateful to the anonymous referee who helped them to improve the first version of this work.

## REFERENCES

[1] ARTZNER, P., DELBAEN, F., EBER, J.-M. and HEATH, D. (1999). Coherent measures of risk. *Math. Finance* **9** 203–228. MR1850791

[2] AVELLANEDA, M., LEVY, A. and PARAS, A. (1995). Pricing and hedging derivative securities in markets with uncertain volatilities. *Appl. Math. Finance* **2** 73–88.

[3] BERBERIAN, S. K. (1974). *Lectures in Functional Analysis and Operator Theory.* Springer, New York. MR0417727

[4] BOULEAU, N. and HIRSCH, F. (1991). *Dirichlet Forms and Analysis on Wiener Space.* de Gruyter, Berlin. MR1133391

[5] CHOQUET, G. (1955). Theory of capacities. *Ann. Inst. Fourier* **5** 131–295. MR0080760

[6] DELBAEN, F. (1992). Representing martingale measures when asset prices are continuous and bounded. *Math. Finance* **2** 107–130.

[7] DELBAEN, F. (2002). Coherent measure of risk on general probability space. In *Advances in Finance and Stochastics* (K. Sandmann and P. J. Schönbucher, eds.) 1–37. Springer, Berlin.

[8] DELLACHERIE, C. (1972). *Capacités et Processus Stochastiques.* Springer, Berlin. MR0448504

[9] DENIS, L. (2005). Solutions of SDE's in a non-dominated model. Preprint, Univ. Maine.

[10] DEPARIS, S. and MARTINI, C. (2004). Superhedging strategies and balayage in discrete time. In *Proceedings of the 4th Ascona Conference on Stochastic Analysis. Random Fields and Applications* **4** 207–222. Birkhäuser, Basel. MR2096290

[11] EL KAROUI, N. and QUENEZ, M.-C. (1995). Dynamic programming and pricing of contingent claims in an incomplete market. *SIAM Control Optim.* **33** 29–66. MR1311659

[12] FEYEL, D. and DE LA PRADELLE, A. (1989). Espaces de Sobolev Gaussiens. *Ann. Inst. Fourier* **39** 875–908. MR1036336

[13] FÖLLMER, H. and KABANOV, YU. M. (1998). Optional decomposition and Lagrange multipliers. *Finance and Stochastics* **2** 69–81. MR1804665

[14] FÖLLMER, H. and SCHIED, A. (2002). Convex measures of risk and trading constraints. *Finance and Stochastics* **6** 429–447. MR1932379

[15] GARSIA, A., RADEMICH, E. and RUMSEY, H. (1970/1971). A real variable lemma and the continuity of paths of some Gaussian processes. *Indiana Univ. Math. J.* **20** 565–578. MR0267632

[16] JACOD, J. and SHIRYAEV, A. (2002). *Limit Theorem for Stochastic Processes.* Springer, Berlin. MR0959133

[17] KRAMKOV, D. (1996). Optional decomposition. *Probab. Theory Related Fields* **105** 459–479. MR1402653

[18] LYONS, T. J. (1995). Uncertain volatility and the risk-free synthesis of derivatives. *Appl. Math. Finance* **2** 117–133.

[19] REVUZ, D. and YOR, M. (1994). *Continuous Martingale and Brownian Motion Integrant,* 2nd ed. Springer, Berlin. MR1303781

Département de Mathématiques
Equipe Analyse et Probabilités
Université d'Evry-Val-d'Essone
Boulevard F. Mitterrand
91025 EVRY Cedex
France
E-mail: ldenis@univ-evry.fr

INRIA and Zeliade Systems
Domaine de Voluceau
Rocquencourt
BP 105
78153 Le Chesnay Cedex
France
E-mail: claude.martini@inria.fr